\documentclass[twoside,a4paper,11pt]{amsart}
\usepackage{epsf,psfig,amscd,amsmath,amssymb,amsxtra,sect}
\usepackage{hyperref}
\def\({\left (}
\def\){\right )}
\newcommand{\GG}{\Gamma}
\newcommand{\GD}{\Delta}
\newcommand{\Gs}{\sigma}
\newcommand{\Gd}{\delta}
\newcommand{\Ge}{\varepsilon}
\newcommand{\R}{\mathbb R}

\newcommand{\sminus}{\smallsetminus}
\newcommand{\p}{\partial}
\newcommand{\Int}{\operatorname{Int}}
\newcommand{\rpp}{\R P^2}
\usepackage[xdvi]{graphicx}

\begin{document}
\title[Dequantization of  Real Algebraic Geometry]
{Dequantization of  Real Algebraic Geometry \\ 
on Logarithmic Paper }

\author{Oleg Viro}
\dedicatory{Uppsala University, Uppsala, Sweden\break
POMI, St.\ Petersburg, Russia}
\address{\newline Department of Mathematics, Uppsala University, 
Box 480,\; S-751 06 Uppsala, Sweden
}
\email{oleg@math.uu.se}

\begin{abstract}
On logarithmic paper some real algebraic curves look like smoothed
broken lines. Moreover, the broken lines can be obtained as limits of
those curves. The corresponding deformation can be  viewed as a
quantization, in which the broken line  is a classical object and the
curves are quantum.  This generalizes to a new connection between
algebraic geometry and the geometry of polyhedra, which is more
straight-forward than the other known connections and gives a new
insight into constructions used in the topology of real algebraic
varieties.   \end{abstract}

\maketitle

\section{Graphs of polynomials on logarithmic
paper}\label{s1}

\subsection{How to visualize a real polynomial?}\label{s1.1} 
If you ever tried to draw the graph for a polynomial of degree greater
than, say, 4 and consisting of at least 4 monomials, you are aware
of the natural difficulties. The graph is too steep. Whatever scale
you choose, either some important details do not fit into the picture or
are too small. The usual recipes from Calculus do not address the
problem, but suggest, instead, to  find roots of the first two
derivatives, which does not seem to be much easier than the original
problem.

\subsection{Logarithmic paper}\label{s1.2}
A physicist or engineer can give a more practical advice, based on
their experience: {\it use (double) logarithmic paper.} 
This is a graph paper, called also {\em log paper}, with a non-uniform net
of coordinate lines and logarithmic scales on both axes. 
On a log paper a point with coordinates
$x$, $y$ is shown at the position with the usual, Cartesian coordinates
equal to $\ln x$, $\ln y$. In other words, the transition to the log 
paper corresponds to the change of coordinates:
$$\begin{cases} u=\ln x\\ v=\ln y.
\end{cases}$$ 
On a log paper the first quadrant is expanded homeomorphically 
to the whole plane, the line $x=1$ occupies the position of the axis of
ordinates, the line $y=1$ occupies the position of the axis of
abscissas, the unit square bounded by these lines and the coordinate
axes occupies the whole third quadrant.  

\subsection{A monomial on logarithmic paper}\label{s1.3}
Let us try to follow the advice to use a log paper. Consider first the 
simplest special case: the graph of a monomial $ax^k$ (i.\ e., the 
curve defined by $y=ax^k$). 
We are forced to consider only positive $x,y$ and hence assume 
$a$ to be positive as well.
Then
$
v=\ln y=\ln(ax^k)=k\ln x+\ln a=ku+\ln a=ku+b
$, where we denote $\ln a$ by $b$. 
Now the curve which we want to draw is defined
in the coordinates $u$ and $v$ by the equation $v=ku+b$.
Everybody knows that this is the straight line with slope
$k$ meeting the axis of ordinates (i.e., $v$-axis) at $(0,b)$.

\subsection{Few slightly more complicated polynomials}\label{s1.4}
\begin{figure}[bhtp] 
\centering
\begin{minipage}[c]{0.34\textwidth}
\centering
\includegraphics[width=1.8in,clip]{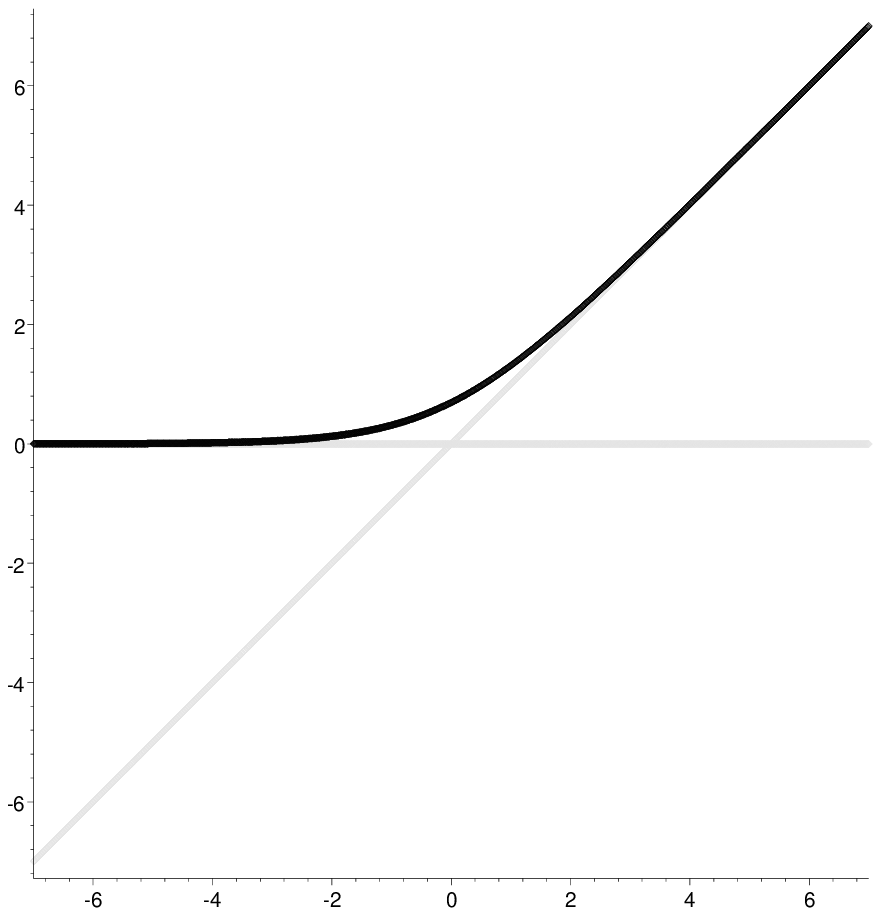}
\end{minipage}%
\begin{minipage}[c]{0.33\textwidth}
\centering
\includegraphics[width=1.65in,clip]{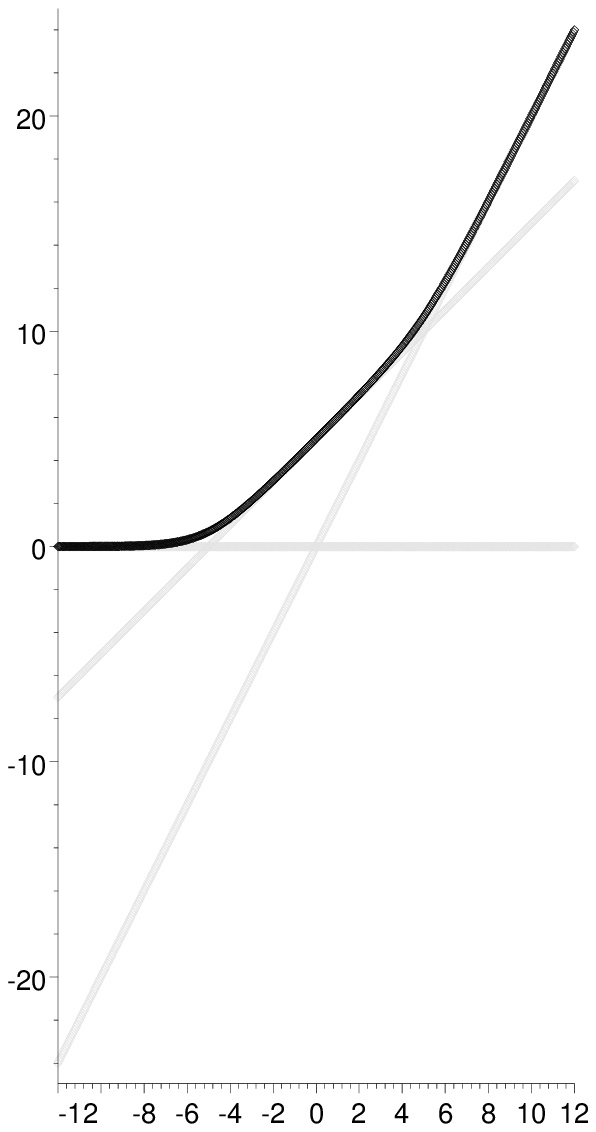}
\end{minipage}%
\begin{minipage}[c]{0.33\textwidth}
\centering
\includegraphics[width=1.65in,clip]{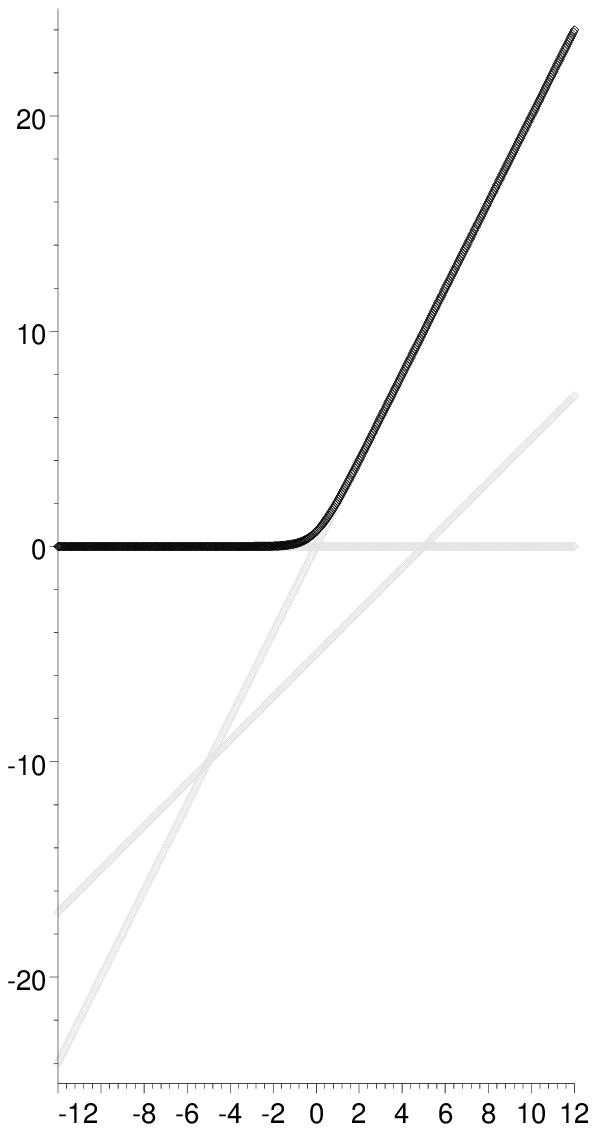}
\end{minipage}
\caption{Log paper graphs of 
$1+x$, $1+e^{5}x+x^2$ and $1+e^{-5}x+x^2$. }
\label{f1}
\end{figure}

First, consider a line $y=1+x$. Then 
$$v=\ln y=\ln(1+x)=\ln(1+e^u).$$
See the left plot in Figure \ref{f1}, where the graph of
$v=\ln(1+e^u)$ is shown together with lines $v=0$ and $v=u$, which
represent on the log paper the monomials $1$ and $x$ involved
in our polynomial $1+x$. The graph of $v=\ln(1+e^u)$ looks like the
broken line $v=\max\{0,u\}$ with  a smoothed corner: it goes along
and above of this broken line getting very close to it as $|u|$ grows.
For $|u|>4$ the difference becomes beyond the resolution.

Here are two further examples: the quadratic polynomials  $1+e^{\pm5}x+x^2$. 
Then 
$$v=\ln(1+e^{\pm5}x+x^2)=\ln(1+e^{u\pm5}+e^{2u}).$$
See the central and right plots in Figure \ref{f1}.
The graph of $v=\ln(1+e^{u\pm5}+e^{2u})$ looks like the broken line 
$v=\max\{0,u\pm5,2u\}$ with smoothed corners. It goes along and above
of this broken line getting very close to it far from its corners.
Notice that the lines $v=0$, $v=u\pm5$ and $v=2u$ represent on the
logarithmic paper the  monomials $1$, $e^{\pm5}x$  and $x^2$, respectively.

\subsection{A polynomial versus the maximum of its
monomials}\label{s1.5}
This suggests, for a polynomial $p(x)=a_nx^{n}
+a_{n-1}x^{n-1}+\dots+a_0$ with positive real coefficients
$a_i=e^{b_i}$, to compare 
the graphs on log paper for $p$ and the maximum 
$M(p)(x)=\max\{a_nx^n,a_{n-1}x^{n-1},\dots,a_0\}$ of its monomials.
Denote the graph on log paper of a function $f$ by $\GG_f$.
With respect to the usual Cartesian coordinates,  $\GG_p$ is the graph of 
$$L_p(u)=\ln\(e^{nu+b_n}+e^{(n-1)u+b_{n-1}}+\dots+e^{b_0}\)$$
and $\GG_{M(p)}$ is the graph of a  piecewise linear convex function
$$
M_p(u)=\max\left\{nu+b_n, \ (n-1)u+b_{n-1},\   \dots\ ,\  b_0  \right\}.
$$

Obviously, $M_p(u)\le L_p(u)\le M_p(u)+\ln(n+1)$. Hence $\GG_{p}$ is 
above $\GG_{M(p)}$, but below a copy of $\GG_{M(p)}$ shifted upwards by 
$\ln(n+1)$. The latter is in fact a rough estimate. It turns to 
equality only at $u$, where all linear functions, whose maximum 
is $M_p(u)$, are equal: $nu+b_n= (n-1)u+b_{n-1}=\dots=b_0$. 

For a generic value of $u$, only one of these functions is equal to
$M_p(u)$. Say $M_p(u)=ku+b_k$, while $M_p(u)>d+lu+b_l$ for some
positive $d$ and each $l\ne k$. Then 
$$L_p(u)<M_p(u)+\ln(1+ne^{-d})<M_p(u)+e^{-d}n.$$

If for some value of $u$ the values of all of the functions $ku+b_k$ 
except two are smaller than $M_p(u)-d$, then 
$$L_p(u)<M_p(u)+\ln(2+(n-1)e^{-d})<M_p(u)+\ln2 +e^{-d}(n-1)/2.$$

Thus, on a logarithmic paper the graph of a generic polynomial with 
positive coefficients lies in a narrow strip along the brocken line
which is the graph of the maximum of its monomials. The width of the
strip is estimated by characteristics of the mutual position of the
lines which are the graphs of the monomials. The less congested the
configuration of these lines, the norrower this strip.

\subsection{Rescalings pushing the graph of a polynomial to a
PL-graph}\label{s1.6}

A natural way to make a configuration of lines less congested 
without changing its topology is to apply a dilation 
$(u,v)\mapsto(Cu,Cv)$ with a large $C>0$. In what follows it is 
more convenient to use instead of $C$ a parameter $h$ related to $C$ 
by $h=1/C$. In terms of $h$
the dilation acts by $(u,v)\mapsto (u/h,v/h)$. It maps the 
graph of $v=ku+b$ to the graph of $v=ku+b/h$. The parallel operation on
monomials replaces $ax^{k}$ by $a^{1/h}x^{k}$. 

Consider the corresponding family of polynomials:
$p_{h}(x)=\sum_{k}a_k^{1/h}x^k$.  On 
log paper, the graphs of its monomials are
obtained  by dilation with ratio $1/h$ from the graphs of the
corresponding 
monomials of $p$. Hence $\GG_{M(p_{h})}$ is the image of $\GG_{M(p)}$
under the same dilation. However, $\GG_{p_{h}}$ is not the image of
$\GG_p$. It still lies in a strip along $\GG_{M(p_{h})}$ and the strip 
is getting narrower as $h$ decreases, but at the corners of 
$\GG_{M(p_{h})}$ the width of the strip cannot become smaller than
$\ln2$.

To keep the picture of our expanding configuration of lines (the graphs
of monomials) independent on $h$, let us make an additional calibration 
of coordinates: set $u_h=hu=h\ln x$,\quad $v_h=hv=h\ln y$. 
Denote by $\GG^h_f$ the graph of a function $y=f(x)$ in the plane with 
coordinates $u_h, v_h$. 

Then $\GG^h_{M(p_{h})}$ does not depend on $h$. The additional 
scaling reduces the width of the strip along  $\GG^h_{M(p_h)}$, 
where $\GG^h_{p_h}$ lies, forcing the width to tend to 0 as $h\to0$.  
Thus $\GG^h_{p_{h}}$ tends to $\GG^h_{M(p_h)}$ (in the $C^0$ sense) as 
$h\to0$.

\section{Quantization}\label{s2}

\subsection{Maslov dequantization of positive real numbers}\label{s2.1}

The rescaling formulas  $u_h=h\ln x$, $v_h=h\ln y$ bring to mind formulas 
related to the Maslov dequantization of real numbers, see e.g.
\cite{LM}, \cite{LMS}. 
The core of the Maslov dequantization is a family
of semirings $\left\{S_h\right\}_{h\in[0,\infty)}$ (recall that a
semiring is a sort of ring, but without subtraction).
As a set, each of $S_h$ is $\R$.
The semiring operations $\oplus_h$ and $\odot_h$ in $S_h$ 
are defined as follows:
\begin{align}a\oplus_h b&=\begin{cases} h\ln(e^{a/h}+e^{b/h}),& \text{ if }h>0 \\
                            \max\{a,b\}, & \text{ if } h=0
\end{cases}\label{oplus}\\
a\odot_h b&= a+b\label{odot}
\end{align}
These operations  depend continuously on $h$.
For each $h>0$ the map 
$$D_h:\R_+\sminus\{0\} \to S_h: x\mapsto h\ln x$$ 
is  a semiring isomorphism of $\left\{\R_+\sminus\{0\},+,\cdot\right\}$ onto
$\left\{S_h,\oplus_h,\odot_h\right\}$, that is
$$D_h(a+b)=D_h(a)\oplus_hD_h(b),\qquad
D_h(ab)=D_h(a)\odot_hD_h(b).
$$
Thus $S_h$ with $h>0$ can be considered as a copy of $\R_+\sminus\{0\}$
with the usual operations of addition and multiplication. On the other
hand, $S_0$ is a copy of $\R$ where the operation of taking maximum is
considered as an addition, and the usual addition, as a multiplication.

Applying the terminology of quantization to this deformation, we must
call $S_0$ a classical object, and $S_h$ with $h\ne0$, quantum ones.
The analogy with Quantum Mechanics motivated  the following {\em
correspondence principle}  formulated by Litvinov and Maslov \cite{LM}
as follows:

{\it  ``There exists a (heuristic) correspondence, in the spirit of the
correspondence principle in Quantum Mechanics, between important,
useful and interesting constructions and results over the field of real
(or complex) numbers  (or the semiring of all nonnegative numbers) and 
similar constructions and results over idempotent semirings.''}

This principle proved to be very fruitful in a number of situations, see
\cite{LM}, \cite{LMS}.
According to the correspondence principle, the idempotent counterpart 
of a polynomial 
$p(x)=a_nx^{n}+a_{n-1}x^{n-1}+\dots+a_0$ is a convex PL-function 
$M_p(u)=\max\left\{nu+b_n, \ (n-1)u+b_{n-1},\ 
\dots\ ,\ b_0\right\}$. As we have seen above, $p$ and $M_p$
are related not only on an heuristic level. In Section \ref{s1.6}
we connected the graph $\GG_p$ of $p$ on logarithmic paper and the
graph $\GG_{M(p)}$ of $M_p$ by a continuous family of graphs
$\{\GG^h_{p_h}\}_{h\in(0,1)}$. 

\subsection{Logarithmic paper as a graphical device for the Maslov
dequantization}\label{s2.2}
As we saw in Section \ref{s1.5}, the graph of a polynomial 
$p(x)=\sum_ka_kx^k$ with positive real coefficients $a_k=e^{b_k}$ 
on log paper is the graph of function $\R\to\R$ defined by 
$v=\ln\left(\sum_ke^{ku+b_k}\right)$. Observe that
$\ln\left(\sum_ke^{ku+b_k}\right)$
is the value in $S_1$ of the polynomial $\sum_kb_kx^k$ at $x=u$. 
Therefore we can identify the graph $\GG_p$ of $p(x)=\sum_ka_kx^k$ on log 
paper with the (Cartesian) graph of the polynomial $\sum_kb_kx^k$ on
$S_1^2$. 

Furthermore, $\GG^h_{p_h}$ is the graph of the function 
$\R\to\R$ defined by 
$$v=h\ln\left(\sum_ka_k^{1/h}e^{(ku)/h}\right)
=h\ln\left(\sum_ke^{(ku+b_k)/h} \right).$$
Observe, that the right hand side is the value in $S_h$ of
the same polynomial $\sum_k b_kx^k$ at $u$.
Therefore we can identify the graph $\GG^h_{p_h}$ of 
$p_h(x)=\sum_ka_k^{1/h}x^k$ on log 
paper with the (Cartesian) graph of $\sum_kb_kx^k$ 
on $S_h^2$.

At last, the graph of $\sum_k b_kx^k$ on $S_0^2$ is the the graph
of $M_p$. 

We see that the whole job of deforming $\GG_p$ to the graph of a 
piecewise linear convex function can be done by the Maslov dequantization:
the deformation consists of the graphs of the same polynomial
$\sum_kb_kx^k$ on $S_h^2$ for $h\in[0,1]$. The coefficients $b_k$ of this
polynomial are logarithms of the coefficients of  the original
polynomial: $b_k=\ln a_k$. Since the map $x\mapsto\ln x:\R_+\sminus0\to
S_1$ was denoted above by $D_1$, we denote by $D_1 F$ the polynomial
obtained from  a polynomial $F$ with positive coefficients by
replacing its coefficients with their logarithms. Thus
$\sum_kb_kx^k=D_1p(x)$. Since $D_1$ is a semiring homomorphism, the
graph $\GG_p$ of $p$ on   log paper is the graph of $D_1p$ on $S_1^2$.
The other graphs involved into the deformation are the graphs of the
same polynomial $D_1p$ on $S_h^2$. They coincide with the graphs on  log
paper of the preimages $p_h$ of $D_1p$ under $D_h$. Indeed,
$p_h(x)=\sum_ka_k^{1/h}x^k$ and
$D_h^{-1}(b_k)=D_h^{-1}D_1(a_k)=e^{D_1(a_k)/h}=e^{(\ln
a_k)/h}=a_k^{1/h}$.

For a real polynomial $p(x)=\sum_ka_kx^k$ with positive coefficients,
we shall call  $p_h(x)=\sum_ka^{1/h}_kx^k$ with $h>0$  the {\em
dequantizing  family\/} of polynomials.

\subsection{Real algebraic geometry as quantized
PL-geometry}\label{s2.4}

The notion of polynomial is central in algebraic geometry. (I 
believe the subject of algebraic geometry would be better described by
the name of {\em polynomial geometry}.) Since a polynomial over $\R$ is
presented so explicitly as a quantization of a piecewise linear convex 
function, one may expect to find along this line explicit relations
between other objects and phenomena of algebraic geometry over $\R$ and
piecewise linear geometry. Indeed, in piecewise linear geometry the
notion of piecewise linear convex function plays almost the same
r\^ole as the notion of polynomial in algebraic geometry.

A representation of real algebraic geometry as a quantized PL-geometry
may be rewarding in many ways. For example, in any quantization there
are {\em classical objects}, i.e., objects which do not change much
under the quantization. Objects of PL-geometry are easier to construct.
If we knew conditions under which a PL object gives rise to a real
algebraic object, which is classical with respect to the Maslov
quantization, then we would have a simple way to construct real
algebraic objects with controlled properties.

\section{Algebraic geometry on logarithmic paper}\label{s3}

\subsection{Speak of real algebraic geometry only positively}\label{s3.1}
Above in Section \ref{s1} we discussed the drawing of graphs on a logarithmic 
paper only for a polynomial with positive coefficients. The graphs 
allowed us to see the behavior of the polynomials only at positive 
values of the argument. This was for a good reason: we used logarithms 
of coordinates. The Maslov dequantization deals only with positive
numbers. Therefore each fragment of algebraic geometry that we want to
dequantize must be reformulated first only in terms of positive numbers.
This seems to be possible for everything belonging to  real algebraic
geometry.

\subsection{Visualizing roots of a polynomial on logarithmic 
paper}\label{s3.2}
Above we could not encounter roots of polynomials, for  a polynomial with
{\em positive\/} coefficients has no {\em positive\/} roots. However if 
we really want to do algebraic geometry on log paper, we must figure
out how to use graphs on log paper for visualizing roots (well,
only {\em positive\/} roots) of an arbitrary real polynomial. 

Any real polynomial $p(x)$ is a difference $p^+(x)-p^-(x)$ of polynomials 
with positive coefficients. We can reformulate the problem of finding 
the positive roots of $p$ as the problem of finding positive values of
$x$ at which $p^+(x)=p^-(x)$. The graphs of $p^+$ and
$p^-$ can be drawn on a log paper, where they are localized in
the strips along broken lines, see Section \ref{s1.5} above. For some
polynomials this picture gives a decent information on the number and
position of the positive roots.

The negative roots of $p(x)$ can be treated in the same way, since
their absolute values are the positive 
roots of $p(-x)$.

\subsection{Plane algebraic curves on logarithmic paper}\label{s3.3}
Now consider a real polynomial $p(x,y)=\sum_{k,l}a_{k,l}x^ky^l$ 
in two variables. Similarly to the case of polynomials in one variable,
in the {\em logarithmic space\/} the graph  of a monomial 
$ax^ky^l$ with
$a>0$ is a plane $w=ku+lv+\ln a$, and the graph $\GG_p$
of a polynomial $p(x,y)$ with positive coefficients lies in a neighborhood 
of a convex piecewise linear surface, which is the graph $\GG_{M(p)}$ of  
the maximum $M(p)$ of the monomials. Furthermore, $p$ is included into
a dequantizing family $p_h$ defined as $\sum_{k,l}a_{k,l}^{1/h}x^ky^l$
for $h>0$,  cf. Section \ref{s2.2}.
The graph
$\GG^h_{p_h}$ of $p_h$ in the logarithmic space with scaled coordinates 
 $u_h=h\ln x$, $v_h=h\ln y$, $w_h=h\ln z$ 
coincides with the graph of the polynomial 
$D_1 p(x,y)=\sum_{k,l}\left(\ln a_{k,l}\right)x^ky^l$ in $S_h^3$. 
These graphs with $h\in(0,1]$ constitute
a continuous deformation of $\GG_p=\GG^1_{p_1}$ to $\GG_{M(p)}$.

For a polynomial $p$ in two variables with arbitrary real coefficients,
denote by $p^+$ the sum of its monomials with positive coefficients,
and put $p^-=p^+-p$. Thus $p$ is canonically presented as a difference
$p^+-p^-$ of two polynomials with positive coefficients. To obtain the
curve  defined on logarithmic paper by the equation $p(x,y)=0$, one
can construct the graphs $\GG_{p^+}$ and $\GG_{p^-}$  for $p^+$ and
$p^-$ in the logarithmic space, which are the surfaces defined
in the usual Cartesian coordinates by  $w=\ln\(p^{\pm}(e^u,e^v)\)$, 
and project the intersection $\GG_{p^+}\cap\GG_{p^-}$ to the plane of
arguments.

For the first approximation of this curve, one may take the broken
line, which is the projection of the intersection of the piecewise
linear surfaces $\GG_{M(p^+)}$ and $\GG_{M(p^-)}$ corresponding to
$p^+$ and $p^-$. 

Of course, it may well happen that the broken line does not even
resemble the curve. This happens to first approximations. However, it 
is very appealing to figure out circumstances under which the broken 
line is a good approximation, for a broken line seems to be much
easier to deal with than an algebraic curve.

\subsection{Constructing algebraic curves, which are classical from our 
quantum point of view}\label{s3.4}

Recall that in the logarithmic space the graph of $ax^ky^l$ 
is a plane
$w=ku+lv+\ln a$. It has a normal vector $(k,l,-1)$ and intersects the 
vertical axis at $(0,0,\ln a)$. Thus if we want to construct a
curve of a given degree $m$, we have to arrange planes whose normals
are fixed: they are $(k,l,-1)$ with integers $k,l$, satisfying
inequalities $0\le k,l,k+l\le m$. The only freedom is in moving them
up and down.  

Consider the pieces of these planes which do not lie under the others. 
They form a convex piecewise linear surface $U$, the graph of the
maximum  of the linear forms defining our planes.  The combinatorial 
structure of faces in $U$ depends on the arrangement. Assume that at
each vertex  of $U$ exactly three of the planes meet. This is a
genericity condition, which can be satisfied by small shifts of the
planes.

Divide now the faces of $U$ arbitrarily into two classes. Denote the
union of one of them by $U^+$, the union of the other by $U^-$.
By genericity of the configuration, the common boundary of $U^+$ and 
$U^-$ is union of disjoint polygonal simple closed curves. It can be
easily realized as the intersection of PL-surfaces $\Gamma_{M(p^+)}$
and $\Gamma_{M(p^-)}$ as above:  take for $p^{\Ge}$ with $\Ge=\pm$ 
the sum of monomials corresponding to the planes of faces forming 
$U^{\Ge}$ and put $p=p^+-p^-$.

Consider now for $1\ge h\ge0$ the curve $C_h\subset S_h^3$ which is the
intersection of the graphs in $S_h^3$ of the polynomials $D_1 p^+$ 
and $D_1p^-$. At $h=0$ this is the intersection  of 
the convex PL-surfaces  $\GG_{M(p^+)}$,  
$\GG_{M(p^-)}$. Due to the genericity condition above, this intersection
is as transversal as one could wish: at all but a finite number of points the
interior part of a face of one of them meets the interior part of a face
of the other one, and at the rest of the points an edge of one of the
surfaces intersects transversaly the interior of a face of the other
surface. 

When $h$ gets positive, the graphs are smoothed, their corners are rounded 
off. The same happens to their intersection curve. While the graphs are
transveral, the intersection curve is deformed isotopically. 

Take the curve corresponding to a value of $h$ such that the transversality
is preserved between 0 and this value. The projection to $(u,v)$-plane
of $C_h$ represents an algebraic curve of degree $m$ on the
scaled logarithmic paper and it can be obtained by a small isotopy 
of the projection of $\p U^+$ to the $(u,v)$-plane.

\subsection{Is this a patchworking?}\label{s3.5}

A construction, which looks similar, has been known in the topology of
real algebraic varieties for about 20 years as  {\em patchworking,} or
{\em Viro's method}.  It has been used to construct
real algebraic varieties with controlled topology and helped to solve a
number of problems. For example, to classify up to isotopy non-singular
real plane projective curves of degree 7 \cite{V1}, \cite{V2} and
disprove the Ragsdale conjecture \cite{It} on the topology of plane curves
formulated \cite{R} as early as  in 1906.
To the best of my knowledge, the patchworking has never been
related to the Maslov quantization.

\section{Patchworking real algebraic curves}\label{s4}

\subsection{The simplest patchworking}\label{s4.1}
Here is a description of a simplified version of patchworking. The
simplifications are of the following 3 kinds:\begin{itemize} 
\item we restrict to the case of nonsingular planar curves, 
\item we assume that all patches are trinomials, and
\item we consider only the part of the curve contained in the
first quadrant (what happens in other quadrants is described soon
after). 
\end{itemize}

{\bf Initial data.}
Let $m$ be a positive integer (it will be the degree of the curve 
under construction) and $\GD$  be the triangle in $\R^2$ with vertices
$(0,0)$, $(m,0)$, $(0,m)$. Let $\tau$ be a {\it convex} triangulation 
of $\GD$ with vertices having integer coordinates. The convexity of $\tau$ 
means that there exists a convex piecewise linear function $\nu:\GD
\longrightarrow {\R_+}$ which is linear on each triangle of $\tau$ and
is not linear on the union of any two triangles of $\tau$. Let the
vertices of $\tau$ be equipped with signs.  The sign (plus or minus) at
the vertex with coordinates $(k,l)$ is denoted by $\Gs_{k,l}$.

{\bf Construction of the piecewise linear curve.}
If a triangle of the triangulation $\tau$ has vertices of different
signs, draw a midline separating pluses from minuses. Denote by $L$ the
union of these midlines. It is a collection of polygonal lines
contained in  $\GD$. The pair $(\GD, L)$ is called the {\it result of
combinatorial patchworking.\/} 

{\bf Construction of polynomials.}
Given initial data $m$, $\GD$, $\tau$
and $\Gs_{k,l}$ as above and a positive
convex function $\nu$ certifying, as above, that the triangulation
$\tau$ is convex.  Consider a one-parameter family of polynomials
\begin{equation}\label{bt}b_t(x,y)=\sum_{\text{\scriptsize$\begin{aligned}&{ 
\text{$(k,l)$  runs over}}\\ &\text{\;  vertices of
$\tau$}\end{aligned}$}}\Gs_{k,l}t^{\nu(k,l)}x^ky^l.\end{equation}
The polynomials $b_t$ are called the results of {\it polynomial
patchworking.\/}
\vskip.6\baselineskip

\noindent
{\bf Patchwork Theorem.} {\it Let
$m$, $\GD$, $\tau$, $\Gs_{k,l}$ and $\nu$ be initial data as above.
Denote by $b_t$ the polynomials obtained by the polynomial
patchworking of these initial data, and by $L$ the
PL-curve in  $\GD$ obtained from the same 
initial data by combinatorial patchworking.

Then for all sufficiently small $t>0$ the polynomial
$b_t$ defines in the first quadrant $\R^2_{++}=\{(x,y)\in\R^2\mid x,y>0\}$
a curve $a_t$ such that the pair
$(\R^2_{++},\; a_t)$ is homeomorphic to the pair
$(\Int\GD,\;L\cap\Int\GD)$.}

\begin{figure}[htb]
\centerline{\includegraphics[clip]{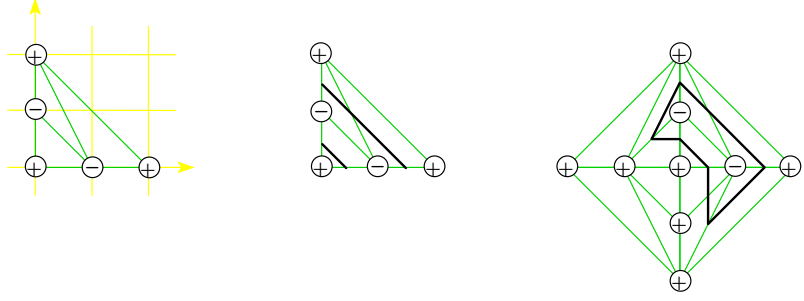}}
\caption{Patchworking: initial data, construction of the PL-curve in the
first quadrant and on the whole plane. The corresponding algebraic
curves are ellipses meeting the coordinate axes in their positive
halves.}
\label{f2}
\end{figure}

\subsection{Patchwork in other quadrants}\label{4.2}

The Patchwork Theorem applied to $b_t(-x,y)$, $b_t(x,-y)$ and $b_t(-x,-y)$
gives a similar topological description of the curve defined in the
other quadrants by $b_t$ with sufficiently small $t>0$. The results can
be collected in the following natural combinatorial construction. 

{\bf Construction of the PL-curve.} Take copies
$\GD_{x} = s_x(\GD)$, $\GD_{y} = s_y(\GD)$, $\GD_{xy} = s(\GD)$
of $\GD$, where $s_x,\; s_y$
are reflections with respect to the coordinate axes and $s= s_x \circ s_y$.
Denote by $A\GD$ the square $\GD \cup \GD_x \cup \GD_y \cup \GD_{xy}$.
Extend the triangulation $\tau$
to a symmetric triangulation of $A\GD$,
and the distribution of signs $\Gs_{i,j}$
to a distribution at the vertices of the extended triangulation by
the following rule: $\Gs_{i,j}\Gs_{\Ge i,\Gd j}\Ge^i\Gd^j=1$, where
$\Ge,\Gd=\pm1$. In other words, passing from a vertex to its
mirror image with respect to an axis we preserve its sign if the
distance from the vertex to the axis is even, and change the sign if
the distance is odd.

If a triangle of the triangulation of $A\GD$ has vertices
of different signs, select (as above) a midline separating pluses from minuses.
Denote by $AL$ the union of the selected
midlines. It is a collection of polygonal lines contained in  $A\GD$.
The pair $(A\GD, AL)$ is called the {\it result of affine combinatorial
patchworking.\/} Glue by $s$ the sides of $A\GD$.  The resulting space
$P\GD$ is homeomorphic to the real projective plane $\rpp$.  Denote
by $PL$ the image of $AL$ in $P\GD$ and call the pair $(P\GD,PL)$ the 
{\it result of projective combinatorial patchworking.\/}
\vskip.6\baselineskip

\noindent
{\bf Addendum to the Patchwork Theorem.} {\it Under the assumptions of Patchwork
Theorem above, for all sufficiently small $t>0$  there exist a
homeomorphism $A\GD\to\R^2$ mapping $AL$ onto the  the affine
curve defined by $b_t$ and a homeomorphism $P\GD\to\R P^2$ mapping $PL$
onto the projective closure of this affine curve.}

\begin{figure}
\centerline{\includegraphics[clip]{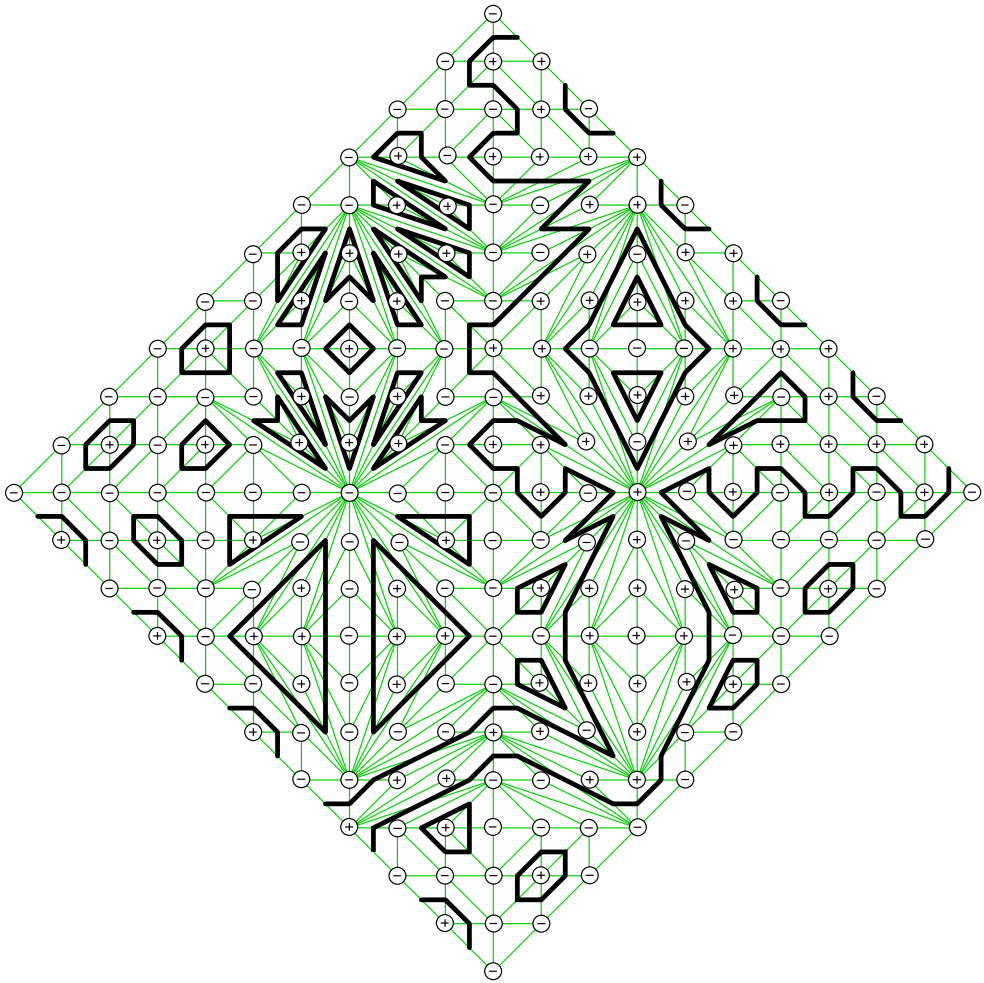}}
\caption{Patchworking of a counter-example to the Ragsdale Conjecture. A curve of
degree 10 with 32 odd ovals constructed by Itenberg \cite{It}.}
\label{f3}
\end{figure}
\subsection{The Simplest Patchworking Coincides With Construction of
Section \ref{s3.4}}\label{s4.3}

The polynomial $b_t$ defined by \eqref{bt} is presented as
$b_t^+-b_t^-$, where 
$$b_t^{\Ge}(x,y)=\sum_{\text{\scriptsize$\begin{aligned}&{ 
\text{$(k,l)$  runs over the vertices}}\\ &\text{\quad of } \tau, \text{ at which $\Gs_{k,l}=\Ge$
}\end{aligned}$}}t^{\nu(k,l)}x^ky^l$$
Observe that polynomials $b_t^{\pm}$ comprise dequantizing families.
Indeed, if we take $p(x,y)=\sum_{k,l}a_{k,l}x^ky^l$ with 
$a_{k,l}=e^{-\nu(k,l)}$, then for $h=-1/{\ln t}$ we obtain
\begin{multline*}
p_h(x,y)=\sum_{k,l}a_{k,l}^{1/h}x^ky^l=\sum_{k,l}e^{-\nu(k,l)/h}x^ky^l=\\
\sum_{k,l}e^{\nu(k,l)\ln t}x^ky^l=\sum_{k,l}t^{\nu(k,l)}x^ky^l.
\end{multline*}
Patchwork Theorem deals with sufficiently small positive $t$, while $h$
in a dequantizing family of polynomials was a small positive number
approaching 0. This is consistent with our setting $h=-1/\ln t$. 

A monomial $a_{k,l}x^ky^l=e^{-\nu(k,l)}x^ky^l$ is presented in the
logarithmic space by the graph of $w=ku+lv-\nu(k,l)$. Hence the graph
of the maximum of linear forms corresponding to all monomials of $p^+$
and $p^-$ is defined by 
\begin{equation}\label{Lt}
w=\max\{ku+lv-\nu(k,l)\;\mid\; (k,l) \text{ runs over vertices
of }\tau \}.
\end{equation}
In \eqref{Lt} we recognize the convex function conjugate to
$\nu$. The graph of \eqref{Lt} is a convex PL surface, whose natural
stratification is dual to the triangulation $\tau$ of $\GD$: the face
which lies on the plane $w=ku+lv-\nu(k,l)$ corresponds to
the vertex $(k,l)$ of $\tau$, two such faces meet at an edge in the
graph of \eqref{Lt} iff the corresponding vertices are
connected with an edge of $\tau$, three faces meet at a vertex 
iff the corresponding vertices of $\tau$ belong to a triangle of
$\tau$. In particular, we see that the configuration of planes
satisfies the genericity condition of Section \ref{s3.4} and 
planes $w=ku+lv-\nu(k,l)$ corresponding to all
monomials of $b_t^{\pm}$ show up in the graph of \eqref{Lt} as its faces.

Some of these faces correspond to monomials of $b_t^+$, the others to
monomials of $b_t^-$. The edges which separate the faces of these two 
kinds constitute a broken line as in Section \ref{s3.4}. 
 These edges are dual to the edges of $\tau$ which
intersect the result $L$ of the combinatorial patchworking. 

Therefore the topology of the projection to the $(u,v)$-plane of 
the broken line coincides with the topology of $L$ in $\GD$.
\qed
\vskip.6\baselineskip
\subsection*{Conclusion}

We see that the quantum point of view (or its graphical log paper
equivalent) gives a natural explanation to
the simplest patchwork construction. The proofs become more conceptual
and straight-forward. Of course, similar but slightly more involved
quantum explanations can be given to all versions of patchwork. 

Let me shortly mention other problems which can be attacked using similar
arguments. 

First of all, this is the Fewnomial Problem. Although
A.\ G.\ Khovansky \cite{Kh} proved that basically all topological characteristics of a real
algebraic variety can be estimated in terms of the number of monomials
in the equations, the known estimates seem to be far weaker than
conjectures. For varieties classical from the quantum point of view a
strong estimates are obvious. It is very compelling to estimate how
much the topology can be complicated by the quantizing deformations.

There seem to be deep relations between the dequantization of algebraic
geometry considered above and the results of I.\ M.\ Gelfand, 
M.\ M.\ Kapranov and A.\ V.\ Zelevinsky on discriminants \cite{GKZ}.
In particular, some monomials in a discriminant are related to 
intersections of hyperplanes in the dequantized polynomial.

Complex algebraic geometry also deserves a dequantization. Especially
relevant may be amoebas introduced in \cite{GKZ}.


\begin{thebibliography}{9999999}
\bibitem{GKZ}  I.\ M.\ Gelfand, M.\ M.\ Kapranov, A.\ V.\ Zelevinsky, 
{\em Discriminants, resultants, and multidimensional determinants}. 
Mathematics: Theory \& Applications. 
Birkh\"auser Boston, Inc., Boston, MA, 1994. x+523 pp. 


\bibitem{It}
I.~Itenberg, Countre-exemples \`a la conjecture de Ragsdale, {\em C.
  R. Acad. Sci. Paris} {\bf 317}, Serie I (1993), 277--282.

\bibitem{Kh}
A.\ G.\ Khovansky, {\em Fewnomials.} Translated from Russian 
by Smilka Zdravkovska. Translations of Mathematical
Monographs, 88. American Mathematical Society, Providence, RI, 1991. 
viii+139 pp. 

\bibitem{LM} Litvinov, G.L. and Maslov, V.P., {\em The correspondence
principle for Idempotent Calculus and some computer applications,\/ }
In book {\it Idempotency} J. Gunawardena (Editor), Cambridge University Press,
Cambridge, 1998, p.420-443.

\bibitem{LMS} G.~L.~Litvinov, V.~P.~Maslov, A.~N.~Sobolevski\u i,
{\em Idempotent Mathematics and Interval Analisys,} Preprint math.SC/9911126,
(1999).

\bibitem{R}
V.~Ragsdale, {\em On the arrangement of the real branches of plane algebraic
  curves}, Amer. J. Math. {\bf 28} (1906), 377--404.


\bibitem{V1}
O.~Ya.~Viro,  Gluing algebraic hypersurfaces and constructions of curves,
  {\em Tezisy Leningradskoj Mezhdunarodnoj Topologicheskoj Konferentsii
1982},
  Nauka, 1983, pp.~149--197 (Russian).

\bibitem{V2}
O.~Ya.~Viro, Gluing of plane real algebraic curves and constructions of
curves of degrees $6$ and $7$, {\em Lecture Notes in Math.}, vol. 1060,
Springer, 1984, pp.~185--200.

\end{thebibliography}
\end{document}